\documentclass[10pt]{amsart}
\usepackage{latexsym}
\usepackage{amssymb}
\usepackage{amsmath}
\usepackage{mathrsfs}
\usepackage[cp866]{inputenc}
\usepackage{bbm}
\usepackage{bbold}
\def\Aff{\operatorname{Aff}}

\def\Fd{\operatorname{Fd}}
\def\spt{\operatorname{supp}}

\begin{document}

\title{
Interaction of Order and Convexity
}

\author{S.~S. Kutateladze}
\address[]{
Sobolev Institute of Mathematics\newline
\indent 4 Koptyug Avenue\newline
\indent Novosibirsk, 630090\newline
\indent Russia}
\email{
sskut@member.ams.org
}
\date{May 28, 2007}

\begin{abstract}
This is an overview  of merging the techniques of
Riesz space theory and convex geometry.
\end{abstract}

\thanks{
Prepared for the
Russian--German geometry meeting dedicated to the 95th anniversary
of A.~D.~Alexandrov (1912--1999)  to take place in St. Petersburg,
June 18--23, 2007.
}
\maketitle
%%%%%%%%%%%%%%%%%%%%%%%%%%%%%%%%%%%%%%%%%%%%%%%%%%%%%%%%%%%%%%%%%%%%%%%%%%

Alexandr Danilovich Alexandrov became
the first and foremost Russian geometer of the twentieth century.
He contributed to mathematics under the slogan:
``Retreat to Euclid,'' remarking that
``the pathos of contemporary mathematics is the return to Ancient Greece.''
Hermann Minkowski revolutionized the theory of numbers
with  the aid of the synthetic geometry of convex surfaces.
The ideas and techniques of the geometry of numbers comprised
the fundamentals of functional analysis which was created by Banach.
The pioneering studies of Alexandrov continued the efforts
of Minkowski and enriched geometry with the methods
of measure theory and functional analysis.
Alexandrov accomplished the turnround to the ancient synthetic geometry
in a much  deeper and subtler sense than it is generally acknowledged today.
Geometry in the large  reduces in no way to overcoming
the local restrictions of differential geometry which bases
upon the infinitesimal methods and ideas of Newton,
Leibniz, and Gauss.

The works of Alexandrov \cite{AD1, AD2} made
tremendous  progress in the theory of
mixed volumes of convex figures. He proved  some fundamental
theorems on convex polyhedra that are celebrated alongside
the theorems of Euler and Cauchy.
While discovering a~solution of the Weyl problem,
Alexandrov suggested  a new synthetic method for proving the theorems of existence.
The results of this research ranked the name of Alexandrov  alongside
the names of Euclid and Cauchy.

Alexandrov enriched the methods of differential geometry by
the tools of functional analysis and measure theory,
driving mathematics to its universal status of the epoch of Euclid.
The mathematics of the ancients was geometry
(there were no other instances of mathematics at all).
Synthesizing geometry with the remaining areas of the today's mathematics,
Alexandrov climbed to the antique ideal of the universal science incarnated in
mathematics. Return to the synthetic methods of {\it mathesis universalis\/}
was inevitable and unavoidable  as well as challenging and fruitful.

\section{Minkowski Duality}

\subsection{}
A~{\it convex figure\/} is a~compact convex set. A~{\it convex body\/}
is a~solid convex figure.
The {\it Minkowski  duality\/} identifies
a~convex figure $S$ in
$\mathbb R^N$ and its {\it support function\/}
$S(z):=\sup\{(x,z)\mid  x\in S\}$ for $z\in \mathbb R^N$.
Considering the members of $\mathbb R^N$ as singletons, we assume that
$\mathbb R^N$ lies in the set $\mathscr V_N$
of all compact convex subsets
of $\mathbb R^N$.
\subsection{}
The classical concept of support function gives rise to abstract convexity
which focuses on the order background of convex sets.

Let $\overline{E}$  be a~complete lattice
 $E$ with the adjoint top $\top:=+\infty$ and bottom $\bot:=-\infty$.
Unless otherwise stated, $Y$ is usually a~{\it Kantorovich space\/}
which is a Dedekind complete vector lattice in another terminology.
Assume further that $H$  is some subset of $E$ which is by implication a~(convex)
cone in $E$, and so the bottom of $E$
lies beyond~$H$. A subset $U$  of~$H$ is {\it convex relative to~}
$H$ or $H$-{\it convex\/}, in symbols $U\in\mathscr{V}(H,\overline{E})$,
provided that $U$ is the $H$-{\it support set\/}
$U^H_p:=\{h\in H\mid h\le p\}$ of some element $p$ of $\overline{E}$.

Alongside the $H$-convex sets we consider
the so-called $H$-convex elements. An element   $p\in \overline{E}$
is  $H$-{\it convex} provided that $p=\sup U^H_p$;~i.e., $p$
represents the supremum of the $H$-support set of~$p$.
The $H$-convex elements comprise the cone which is denoted by
${\mathscr Cnv}(H,\overline{E}$).  We may omit  the references to $H$ when $H$ is clear
from the context. It is worth noting that
convex elements and sets are ``glued together''
by the {\it Minkowski duality\/} $ \varphi:p\mapsto U^H_p$.
This duality enables us to study convex elements and sets simultaneously.

Since the classical results by Fenchel \cite{Fenchel}
and H\"ormander
\cite{Her, Notions} we know
that the most convenient and conventional classes of convex  functions
and sets are
${\mathscr Cnv}(\Aff(X),\overline{\mathbb R^X} )$ and
$\mathscr{V}(X',\overline{\mathbb R^X})$.
Here    $X$ is a locally convex space,   $X'$ is the dual of~$X$,
and $\Aff(X)$ is the space of affine functions on  $X$
(isomorphic with  $X'\times \mathbb R$).

In the first case the Minkowski duality is the mapping
 $f\mapsto\text{epi} (f^*)$ where
$$
f^*(y):=\sup\nolimits_{x\in X}(\langle y,x\rangle - f(x))
$$
is the {\it Young--Fenchel transform\/} of~$f$ or
the {\it conjugate function\/} of~$f$.
In the second case we return to the classical identification
of $U$  in $\mathscr{V}(X',\overline{\mathbb R}^X)$
and the standard  support function that uses
the canonical pairing $\langle\cdot,\cdot\rangle$ of~$X'$ and~$X$.

This idea of abstract convexity lies behind many current objects
of analysis and geometry. Among them we list the ``economical'' sets
with boundary points meeting the Pareto criterion, capacities, monotone
seminorms, various classes of functions convex in some generalized sense,
for instance, the Bauer convexity in Choquet  theory, etc.
It is curious that there are ordered vector spaces consisting of
the convex elements with respect to  narrow cones with finite generators.
Abstract convexity is traced and reflected, for instance, in
\cite{MD}--%
%, KutRub, Singer, PR, Rub,
\cite{Rub}.

\section{Positive Functionals over Convex Objects}

\subsection{}
The Minkowski duality makes $\mathscr V_N$ into a~cone
in the space $C(S_{N-1})$  of continuous functions on the Euclidean unit sphere
$S_{N-1}$, the boundary of the unit ball $\mathfrak z_N$.
This yields  the so-called {\it  Minkowski structure\/} on $\mathscr V_N$.
Addition of the support functions
of convex figures amounts to taking their algebraic sum, also called the
{\it Minkowski addition}. It is worth observing that the
{\it linear span}
$[\mathscr V_N]$ of~$\mathscr V_N$ is dense in $C(S_{N-1})$, bears
a~natural structure of a~vector lattice
and is usually referred to as the {\it space of convex sets}.
The study of this space stems from the pioneering breakthrough of
Alexandrov in 1937 and the further insights of  Radstr\"{o}m \cite{Rad},
H\"{o}rmander \cite{Her}, and Pinsker \cite{Pin}.

\subsection{}
It was long ago in 1954 that  Reshetnyak suggested
in his Ph.~D. thesis \cite{Resh} to compare positive measures
on~$S_{N-1}$ as follows.

 A measure $\mu$ {\it linearly majorizes} or {\it dominates}
a~measure $\nu$ provided that to each decomposition of
$S_{N-1}$ into finitely many disjoint Borel sets $U_1,\dots,U_m$
there are measures $\mu_1,\dots,\mu_m$ with sum $\mu$
such that every difference $\mu_k - \nu|_{U_k}$ annihilates
all restrictions to $S_{N-1}$ of linear functionals over
$\mathbb R^N$. In symbols, we write $\mu\,{\gg}{}_{\mathbb R^N} \nu$.

Reshetnyak    proved that
$$
\int_{S_N-1} p d\mu \ge  \int_{S_N-1} p d\nu
$$
for each  sublinear    functional  $p$
on  $\mathbb R^N$   if   $\mu\,{\gg}{}_{\mathbb R^N} \nu$.
This gave an important trick for generating positive linear functionals
over various classes of convex  surfaces and functions.

\subsection{}
A~similar idea was suggested by   Loomis \cite{Loomis}
in 1962 within Choquet theory:

A~measure $\mu$ {\it affinely majorizes\/} or {\it dominates\/}
a measure $\nu$, both given     on a compact convex subset $Q$ of a locally convex space $X$,
provided that     to each decomposition of
$\nu$ into finitely many summands
$\nu_1,\dots,\nu_m$  there are measures $\mu_1,\dots,\mu_m$
whose sum is $\mu$ and for which every difference
$\mu_k - \nu_k$ annihilates all restrictions
to  $Q$  of affine   functionals over $X$.
In symbols, $\mu\,{\gg}{}_{ \Aff(Q)} \nu$.

Cartier, Fell, and Meyer \cite{CFM} proved in 1964 that
$$
\int_{Q} f d\mu \ge  \int_{Q} f d\nu
$$
for each continuous convex function  $f$
on  $Q$   if and only if   $\mu\,{\gg}{}_{\Aff(Q)} \nu$.
An analogous necessity part for linear majorization was published
in 1970, cf. \cite{Kut70}.

\subsection{} Majorization is a vast subject \cite{Marshall}.
We only site  one of the relevant abstract claims of subdifferential
calculus~\cite{Subdiff}:

\subsection{\bf Theorem}
Assume that $H_1,\dots,H_N$ are cones in a Riesz space~$X$.
Assume further that $f$ and $g$ are positive functionals on~$X$.
The inequality
$$
f(h_1\vee\dots\vee h_N)\ge g(h_1\vee\dots\vee h_N)
$$
holds for all
$h_k\in H_k$ $(k:=1,\dots,N)$
if and only if to each decomposition
of~$g$ into a sum of~$N$ positive terms
$g=g_1+\dots+g_N$
there is a decomposition of~$f$ into a sum of~$N$
positive terms $f=f_1+\dots+f_N$
such that
$$
f_k(h_k)\ge g_k(h_k)\quad
(h_k\in H_k;\ k:=1,\dots,N).
$$

\section{Alexandrov Measures and the Blaschke Structure}
The celebrated {\it Alexandrov Theorem\/} \cite[p.~108]{AD1} proves the unique existence of
a translate of a convex body given its surface area function.
Each surface area function is an {\it Alexandrov measure}.
So we call a positive measure on the unit sphere which is supported by
no great hypersphere and which annihilates
singletons.
The last property of a measure is referred to as translation
invariance in the theory of convex surfaces. Thus,
each Alexandrov measure is a translation-invariant
additive functional over the cone
$\mathscr V_N$.

This yields
some abstract cone structure
that results from identifying the coset of translates
$\{z+\mathfrak x \mid  z\in \mathbb R^N\}$
of a  convex body $\mathfrak x$
the corresponding Alexandrov measure on the unit sphere which we call the
{\it surface area function\/} of the coset of $\mathfrak x$ and
denote by  $\mu (\mathfrak x)$.
The soundness of this parametrization rests on the Alexandrov Theorem.

The cone of positive translation-invariant measures in the
dual $C'(S_{N-1})$ of
 $C(S_{N-1})$ is denoted by~$\mathscr A_N$.
We now agree on some preliminaries.

Given $\mathfrak x, \mathfrak y\in \mathscr V_N$, we let the record
$\mathfrak x\,{=}{}_{\mathbb R^N}\mathfrak y$ mean that $\mathfrak x$
and $\mathfrak y$  are  equal up to translation or, in other words,
are translates of one another.
We may say that ${=}{}_{\mathbb R^N}$ is the associate equivalence of
the preorder $\ge{}_{\mathbb R^N}$ on $\mathscr V_N$  which symbolizes
the possibility of inserting one figure into the other
by translation.
Arrange the factor set $\mathscr V_N/\mathbb R^N$ which consists of
the cosets of translates of the members of $\mathscr V_N$.
Clearly,
$\mathscr V_N/\mathbb R^N$ is a cone in the factor space $[\mathscr V_N]/\mathbb R^N$
of the vector space $[\mathscr V_N]$ by the subspace $\mathbb R^N$.
There is a natural bijection between $\mathscr V_N/\mathbb R^N$ and $\mathscr A_N$.
Namely,  we identify the coset of singletons with the zero measure.
To the straight line segment with endpoints $x$ and $y$,
we assign the measure
$
|x-y|(\varepsilon _{(x-y)/|x-y|}+\varepsilon _{(y-x)/|x-y|}),
$
where $|\,\cdot\,|$ stands for the Euclidean norm
and the symbol
${\varepsilon}_z $
for $z\in S_{N-1}$
stands for the {\it Dirac measure} supported at~$z$.
If the dimension of the affine span
$\Aff(\mathfrak x)$ of a representative $\mathfrak x$ of a coset in
$\mathscr V_N/\mathbb R^N$
is greater than unity, then we assume that $\Aff(\mathfrak x)$ is a subspace
of $\mathbb R^N$
and identify this class with the surface area function
of~$\mathfrak x$ in $\Aff(\mathfrak x)$ which is
some measure on $S_{N-1}\cap\Aff(\mathfrak x)$
in this event.
Extending the measure by zero to a measure on  $S_{N-1}$,
we obtain the member of  $\mathscr A_N$ that we
assign to the coset of all translates of $\mathfrak x$.
The fact that this correspondence is one-to-one follows easily from
the Alexandrov Theorem.

The vector space structure on the set of regular Borel measures
induces
in $\mathscr A_N$ and, hence, in
$\mathscr V_N/\mathbb R^N$
the structure of an abstract
cone or, strictly speaking, the structure
of
a commutative  $\mathbb R_+$-operator
semigroup with cancellation.
This structure on $\mathscr V_N/\mathbb R^N$ is called the
{\it  Blaschke structure} (cp.~\cite{Bla} and the references therein).
Note that the sum of the surface area functions
of $\mathfrak x$ and $\mathfrak y$ generates a unique class
$\mathfrak x\# \mathfrak y$ which is referred to as the
{\it Blaschke sum\/} of $\mathfrak x$ and~$\mathfrak y$.

Let $C(S_{N-1})/\mathbb R^N$ stand for the factor space of
$C(S_{N-1})$ by the subspace of all restrictions of linear
functionals on $\mathbb R^N$ to $S_{N-1}$.
Denote by $[\mathscr A_N]$ the space $\mathscr A_N-\mathscr A_N$
of translation-invariant measures. It is easy to see
that $[\mathscr A_N]$ is also the linear span
of the set of Alexandrov measures.
The spaces $C(S_{N-1})/\mathbb R^N$ and $[\mathscr A_N]$ are made dual
by the canonical bilinear form
$$
\langle f ,\mu\rangle={1\over N}\int\nolimits_{S_{N-1}}fd\mu \quad
(f\in C(S_{N-1})/\mathbb R^N,\ \mu \in[\mathscr A_N]).
$$
For $\mathfrak x\in\mathscr V_N/\mathbb R^N$ and $\mathfrak y\in\mathscr A_N$,
the quantity
$\langle {\mathfrak x} ,{\mathfrak y}\rangle$ coincides with the
{\it mixed volume\/}
$V_1 (\mathfrak y ,\mathfrak x)$.
The space $[\mathscr A_N]$ is usually furnished with
the weak topology induced by the above indicated duality
with $C(S_{N-1})/\mathbb R^N$.

\section{Cones of Feasible Directions}
\subsection{}
By the {\it dual\/} $K^*$ of a given cone $K$
in a vector space $X$ in duality with another vector space
$Y$, we mean the set of all positive linear functionals on
$K$; i.e.,
$K^*:=\{y\in Y\mid (\forall x\in K)\ \langle x ,y\rangle\ge 0\}$.
Recall also that
to a convex subset $U$ of $X$ and a point $\bar x$ in  $U$
there corresponds the cone
$$
U_{\bar x}:=\Fd (U,\bar x):=\{h\in X\mid (\exists \alpha \ge 0)\ \bar x+\alpha h\in U \}
$$
which is called the {\it cone of feasible directions\/}
of $U$ at $\bar x$.
Fortunately, description is available for
all dual cones we need.

\subsection{}
Let $\bar {\mathfrak x}\in{\mathscr A}_N$.
Then the dual  $\mathscr A^*_{N,\bar{\mathfrak x}}$ of the cone of
feasible directions of $\mathscr A_Nn$
at~$\bar{\mathfrak x}$ may be represented as follows
$$
{\mathscr A}^*_{N,\bar{\mathfrak x}}=\{f\in{\mathscr A}^*_N\mid
\langle\bar {\mathfrak x},f\rangle=0\}.
$$

\subsection{}
Let $\mathfrak x$ and $\mathfrak y$ be  convex figures. Then

(1) $\mu(\mathfrak x)- \mu(\mathfrak y)\in \mathscr V^*_N
\leftrightarrow \mu(\mathfrak x)\,{\gg}{}_{\mathbb R^N} \mu(\mathfrak y)$;

(2) If $\mathfrak x\ge{}_{\mathbb R^N}\mathfrak y$
then  $\mu(\mathfrak x)\,{\gg}{}_{\mathbb R^N} \mu(\mathfrak y)$;

(3) $\mathfrak x\ge{}_{\mathbb R^2}\mathfrak y\leftrightarrow
\mu(\mathfrak x)\,{\gg}{}_{\mathbb R^2} \mu(\mathfrak y)$;

(4) If $\mathfrak y-{\bar{\mathfrak x}}\in\mathscr A^*_{N,\bar{\mathfrak x}}$ then
$\mathfrak y=_{\mathbb R^N}\bar{\mathfrak x}$;

(5) If $\mu (\mathfrak y)-\mu (\bar{\mathfrak x})\in\mathscr V^*_{N,\bar{\mathfrak x}}$
then
$\mathfrak y=_{\mathbb R^N}\bar{\mathfrak x}$.

\noindent
It stands to reason to avoid  discriminating between a  convex figure,
the respective coset of translates in  $\mathscr V_N/\mathbb R^N$,
and the corresponding measure in $\mathscr A_N$.

\section{Comparison Between  the Blaschke and Minkowski Structures}

The isoperimetric-type problems with
subsidiary constraints on location of convex figures
comprise  a unique class of meaningful
extremal problems  with two essentially different parametrizations.
The principal features of the latter are seen
from
the table.

$$
\vbox{\tabskip=0pt \offinterlineskip
\def\tablerule{\noalign{\hrule}}
\halign to 325pt{ \vrule#\tabskip=1em plus2em&\strut
\hfil#& \vrule#& \hfil#& \vrule#& \hfil#& \vrule#
\tabskip=0pt\cr\tablerule
height5pt&
\omit &&  \omit && \omit &\cr
&\omit\hidewidth {\scshape Object} \hidewidth&&\omit
\hidewidth {\scshape Minkowski's}  \hidewidth&&
\omit\hidewidth {\scshape Blaschke's} \hidewidth&\cr
height2pt&\omit &&  \omit && \omit &\cr
&\omit\hidewidth {\scshape of Parametrization}\hidewidth&&\omit
\hidewidth{\scshape  Structure} \hidewidth&&
\omit\hidewidth {\scshape Structure}\hidewidth&\cr
height2pt& \omit && \omit && \omit &\cr
\tablerule
height2pt& \omit &&  \omit && \omit &\cr
& cone of sets\hfil &&${\mathscr V}_N/\mathbb R^N$\hfil&&${\mathscr A}_N$\hfil&\cr
height4pt& \omit  &&  \omit && \omit &\cr
&dual cone\hfil &&${\mathscr V}^*_N$\hfil &&${\mathscr A}^*_N$\hfil&\cr
height4pt& \omit  &&  \omit && \omit &\cr
&positive cone\hfil &&$\mathscr A^*_N $\hfil&&$\mathscr A_N$\hfil&\cr
height4pt& \omit  &&  \omit && \omit &\cr
&typical linear\hfil  &&$V_1 (\mathfrak z_N,\,\cdot\,)$\hfil&& $V_1(\,\cdot\,,\mathfrak z_N)$\hfil&\cr
&functional\hfil  && (width)\hfil && (area)\hfil &\cr
height4pt& \omit  &&  \omit && \omit &\cr
&concave functional\hfil  &&$V^{1/N}(\,\cdot\,)$\hfil
&&$V^{(N-1)/N}(\,\cdot\,)$\hfil &\cr
&(power of volume)\hfil  &&  &&  &\cr
height4pt& \omit  &&  \omit && \omit &\cr
&simplest convex\hfil  && isoperimetric\hfil && Urysohn's\hfil &\cr
&program\hfil && problem\hfil && problem\hfil &\cr
height4pt& \omit  &&  \omit && \omit &\cr
&operator-type\hfil  && inclusion\hfil && inequalities\hfil &\cr
&constraint\hfil  && of figures\hfil && on ``curvatures''\hfil &\cr
height4pt& \omit  &&  \omit && \omit &\cr
&Lagrange's multiplier\hfil  &&  surface\hfil && function\hfil &\cr
height4pt& \omit  &&  \omit && \omit &\cr
&differential of volume\hfil  &&  &&  &\cr
&at a point $\bar{\mathfrak x} $ \hfil  && && &\cr
&is proportional to\hfil&&$V_1(\bar{\mathfrak x},\,\cdot\,) $\hfil&& $V_1(\,\cdot\,,\bar{\mathfrak x})$\hfil&\cr
height3pt& \omit  &&  \omit && \omit &\cr
\tablerule
\noalign{\smallskip}\cr }}
$$

This table shows that the classical isoperimetric problem
is not a convex program in the Minkowski structure for
$N\ge 3$. In this event a necessary optimality condition
leads to a solution only under extra regularity conditions.
Whereas in the Blaschke structure this problem is a convex program
whose optimality criterion reads: ``Each solution is a ball.''

The problems are challenging that contain some constrains of inclusion type:
for instance, the isoperimetric problem or Urysohn problem with the requirement that
the solutions lie among the subsets or supersets of a~given body.
These problems can be
solved in a generalized sense, ``modulo'' the Alexandrov Theorem.
These problems can be
solved in a generalized sense ``modulo'' the Alexandrov Theorem.
Clearly, some convex combination  of the ball and a tetrahedron is proportional
to the solution of the  Urysohn problem in this tetrahedron.
If we replace the condition on the integral which
is characteristic of the Urysohn problem \cite{Ury, Lyus}
by a~constraint on the surface area or other mixed volumes
of a more general shape then we come to possibly nonconvex
programs for which a similar reasoning yields
only necessary extremum conditions in general.
Recall that  in case $N=2$ the Blaschke sum transforms as usual into the
Minkowski sum modulo translates.

The task of choosing an appropriate parametrization
for a wide class of problems is practically unstudied
in general. In particular, those problems of geometry
remain unsolved which combine constraints each of which is linear
in one of the two vector structures on the set of  convex figures.
The simplest example of an unsolved  ``combined''
problem is the internal isoperimetric problem
in the space  $\mathbb R^N$ for $N\ge 3$. The only instance of progress
is due to Pogorelov  who found in \cite{Pog} the
form of a soap bubble inside a~three-dimensional tetrahedron.
This happens to be proportional to the Minkowski convex combination of
the ball and the solution to the internal Urysohn problem in the
tetrahedron.

The above geometric facts make it reasonable to
address the general problem of parametrizing
the important classes of extremal problems of practical provenance.

\subsection{}
By way of example, consider the {\it external Urysohn problem}:
Among the convex figures, circumscribing $\mathfrak x_0 $ and having
integral width fixed, find a convex body of greatest volume.

\subsection{\bf Theorem}
A feasible convex body $\bar {\mathfrak x}$ is a solution
to~the external Urysohn problem
if and only if there are a positive  measure~$\mu $
and a positive real $\bar \alpha \in \mathbb R_+$ satisfying

(1) $\bar \alpha \mu (\mathfrak z_N)\,{\gg}{}_{\mathbb R^N}\mu (\bar {\mathfrak x})+\mu $;

(2)~$V(\bar {\mathfrak x})+{1\over N}\int\nolimits_{S_{N-1}}
\bar {\mathfrak x}d\mu =\bar \alpha V_1 (\mathfrak z_N,\bar {\mathfrak x})$;

(3)~$\bar{\mathfrak x}(z)={\mathfrak x}_0 (z)$
for all $z$ in the support of~$\mu $.

\subsection{} If, in particular, ${\mathfrak x}_0 ={\mathfrak z}_{N-1}$ then the sought body
is a {\it spherical lens}, that is, the intersection
of   two balls of the same radius; while
the critical measure is the restriction
of the surface area function
of the ball of radius
$\bar \alpha ^{1/(N-1)}$
to the complement of the support of the lens to~$S_{N-1}$.
If ${\mathfrak x}_0 ={\mathfrak z}_1 $ and $N=3$ then our result implies that
we should seek
a solution  in the class of the so-called
spindle-shaped constant-width surfaces of revolution.

\subsection{}
We turn now to consider the {\it internal Urysohn problem
with a current hyperplane} (cp.~\cite{Kut02}):
Find two convex figures $\bar{\mathfrak x}$ and $\bar{\mathfrak y}$
lying in a given convex body
$\mathfrak x_o$,
 separated by a hyperplane with the unit outer normal $z_0$,
and having the greatest total volume
of $\bar{\mathfrak x}$ and $\bar{\mathfrak y}$
given the sum of their integral widths.

\subsection{\bf Theorem}
A feasible pair of convex bodies $\bar{\mathfrak x}$ and $\bar{\mathfrak y}$
solves the internal Urysohn problem with a current hyperplane
if and only if
there are convex figures  $\mathfrak x$ and $\mathfrak y$
and positive reals
$\bar\alpha $ and $\bar\beta$  satisfying

{\rm(1)}  $\bar{\mathfrak x}=\mathfrak x \# \bar\alpha\mathfrak z_N$;

{\rm(2)} $\bar{\mathfrak y}=\mathfrak y \# \bar\alpha\mathfrak z_N$;

{\rm(3)} $\mu(\mathfrak x) \ge \bar\beta\varepsilon_{z_0} $, $\mu(\mathfrak y) \ge \bar\beta\varepsilon_{-z_0} $;

{\rm(4)} $\bar {\mathfrak x}(z)=\mathfrak x_0 (z)$ for all $z\in \spt(\mathfrak x)\setminus \{z_0\} $;

{\rm(5)} $\bar {\mathfrak y}(z)=\mathfrak x_0 (z)$ for all $z\in \spt(\mathfrak x)\setminus \{-z_0\} $,

\noindent
with $\spt(\mathfrak x)$ standing for the {\it support\/} of $\mathfrak x$,
i.e. the support of the surface area measure $\mu(\mathfrak x)$
of~$\mathfrak x$.

The internal isoperimetric problem and its analogs
seem indispensable since we have no adequate
means for expressing their solutions. The new level of understanding
is in order in convexity that we may hope to achieve
with the heritage of  Alexandrov, the teacher of universal freedom
in geometry.

\bibliographystyle{plain}

\enddocument